\magnification=1200

\hsize=11.25cm    
\vsize=18cm     
\parindent=12pt   \parskip=5pt     

\hoffset=.5cm   
\voffset=.8cm   

\pretolerance=500 \tolerance=1000  \brokenpenalty=5000

\catcode`\@=11

\font\eightrm=cmr8         \font\eighti=cmmi8
\font\eightsy=cmsy8        \font\eightbf=cmbx8
\font\eighttt=cmtt8        \font\eightit=cmti8
\font\eightsl=cmsl8        \font\sixrm=cmr6
\font\sixi=cmmi6           \font\sixsy=cmsy6
\font\sixbf=cmbx6

\font\tengoth=eufm10 
\font\eightgoth=eufm8  
\font\sevengoth=eufm7      
\font\sixgoth=eufm6        \font\fivegoth=eufm5

\skewchar\eighti='177 \skewchar\sixi='177
\skewchar\eightsy='60 \skewchar\sixsy='60

\newfam\gothfam           \newfam\bboardfam

\def\tenpoint{
  \textfont0=\tenrm \scriptfont0=\sevenrm \scriptscriptfont0=\fiverm
  \def\rm{\fam\z@\tenrm}
  \textfont1=\teni  \scriptfont1=\seveni  \scriptscriptfont1=\fivei
  \def\oldstyle{\fam\@ne\teni}\let\old=\oldstyle
  \textfont2=\tensy \scriptfont2=\sevensy \scriptscriptfont2=\fivesy
  \textfont\gothfam=\tengoth \scriptfont\gothfam=\sevengoth
  \scriptscriptfont\gothfam=\fivegoth
  \def\goth{\fam\gothfam\tengoth}
  
  \textfont\itfam=\tenit
  \def\it{\fam\itfam\tenit}
  \textfont\slfam=\tensl
  \def\sl{\fam\slfam\tensl}
  \textfont\bffam=\tenbf \scriptfont\bffam=\sevenbf
  \scriptscriptfont\bffam=\fivebf
  \def\bf{\fam\bffam\tenbf}
  \textfont\ttfam=\tentt
  \def\tt{\fam\ttfam\tentt}
  \abovedisplayskip=12pt plus 3pt minus 9pt
  \belowdisplayskip=\abovedisplayskip
  \abovedisplayshortskip=0pt plus 3pt
  \belowdisplayshortskip=4pt plus 3pt 
  \smallskipamount=3pt plus 1pt minus 1pt
  \medskipamount=6pt plus 2pt minus 2pt
  \bigskipamount=12pt plus 4pt minus 4pt
  \normalbaselineskip=12pt
  \setbox\strutbox=\hbox{\vrule height8.5pt depth3.5pt width0pt}
  \let\bigf@nt=\tenrm       \let\smallf@nt=\sevenrm
  \normalbaselines\rm}

\def\eightpoint{
  \textfont0=\eightrm \scriptfont0=\sixrm \scriptscriptfont0=\fiverm
  \def\rm{\fam\z@\eightrm}
  \textfont1=\eighti  \scriptfont1=\sixi  \scriptscriptfont1=\fivei
  \def\oldstyle{\fam\@ne\eighti}\let\old=\oldstyle
  \textfont2=\eightsy \scriptfont2=\sixsy \scriptscriptfont2=\fivesy
  \textfont\gothfam=\eightgoth \scriptfont\gothfam=\sixgoth
  \scriptscriptfont\gothfam=\fivegoth
  \def\goth{\fam\gothfam\eightgoth}
  
  \textfont\itfam=\eightit
  \def\it{\fam\itfam\eightit}
  \textfont\slfam=\eightsl
  \def\sl{\fam\slfam\eightsl}
  \textfont\bffam=\eightbf \scriptfont\bffam=\sixbf
  \scriptscriptfont\bffam=\fivebf
  \def\bf{\fam\bffam\eightbf}
  \textfont\ttfam=\eighttt
  \def\tt{\fam\ttfam\eighttt}
  \abovedisplayskip=9pt plus 3pt minus 9pt
  \belowdisplayskip=\abovedisplayskip
  \abovedisplayshortskip=0pt plus 3pt
  \belowdisplayshortskip=3pt plus 3pt 
  \smallskipamount=2pt plus 1pt minus 1pt
  \medskipamount=4pt plus 2pt minus 1pt
  \bigskipamount=9pt plus 3pt minus 3pt
  \normalbaselineskip=9pt
  \setbox\strutbox=\hbox{\vrule height7pt depth2pt width0pt}
  \let\bigf@nt=\eightrm     \let\smallf@nt=\sixrm
  \normalbaselines\rm}

\tenpoint

\def\pc#1{\bigf@nt#1\smallf@nt}         \def\pd#1 {{\pc#1} }

\frenchspacing

\def\raggedbottom{\topskip 10pt plus 36pt\r@ggedbottomtrue}

\def\pointir{\unskip . --- \ignorespaces}

\def\Medbreak{\vskip-\lastskip\medbreak}

\long\def\th#1 #2\enonce#3\endth{
   \Medbreak\noindent
   {\pc#1} {#2\unskip}\pointir{\it #3}\smallskip}

\def\decale#1{\smallbreak\hskip 28pt\llap{#1}\kern 5pt}
\def\decaledecale#1{\smallbreak\hskip 34pt\llap{#1}\kern 5pt}
\def\puce{\smallbreak\hskip 6pt{$\scriptstyle\bullet$}\kern 5pt}

\def\eqalign#1{\null\,\vcenter{\openup\jot\m@th\ialign{
\strut\hfil$\displaystyle{##}$&$\displaystyle{{}##}$\hfil
&&\quad\strut\hfil$\displaystyle{##}$&$\displaystyle{{}##}$\hfil
\crcr#1\crcr}}\,}

\catcode`\@=12

\showboxbreadth=-1  \showboxdepth=-1

\newcount\numerodesection \numerodesection=1
\def\section#1{\bigbreak
 {\bf\number\numerodesection.\ \ #1}\nobreak\medskip
 \advance\numerodesection by1}

\mathcode`A="7041 \mathcode`B="7042 \mathcode`C="7043 \mathcode`D="7044
\mathcode`E="7045 \mathcode`F="7046 \mathcode`G="7047 \mathcode`H="7048
\mathcode`I="7049 \mathcode`J="704A \mathcode`K="704B \mathcode`L="704C
\mathcode`M="704D \mathcode`N="704E \mathcode`O="704F \mathcode`P="7050
\mathcode`Q="7051 \mathcode`R="7052 \mathcode`S="7053 \mathcode`T="7054
\mathcode`U="7055 \mathcode`V="7056 \mathcode`W="7057 \mathcode`X="7058
\mathcode`Y="7059 \mathcode`Z="705A


\def\diagram#1{\def\normalbaselines{\baselineskip=0pt\lineskip=5pt}
\matrix{#1}}

\def\ufl#1#2#3{\llap{$\textstyle #1$}
\left\uparrow\vbox to#3{}\right.\rlap{$\textstyle #2$}}

\def\hfl#1#2#3{\smash{\mathop{\hbox to#3{\rightarrowfill}}\limits
^{\textstyle#1}_{\textstyle#2}}}

\def\pgoth{{\goth p}}

\def\Q{{\bf Q}}

\def\N{{\bf N}}

\def\Z{{\bf Z}}

\def\F{{\bf F}}

\def\Hom{\mathop{\rm Hom}\nolimits}

\def\Card{\mathop{\rm Card}\nolimits}
\def\Gal{\mathop{\rm Gal}\nolimits}
\def\Ker{\mathop{\rm Ker}\nolimits}

\def\series#1{(\!(#1)\!)}

\def\to{\rightarrow}

\def\mod{\mathop{\rm mod.}\nolimits}
\def\pmod#1{\;(\mod#1)}

\def\char{\mathop{\rm char}\nolimits}

\def\Aut{\mathop{\rm Aut}\nolimits}

\def\boxit#1{\vbox{\hrule\hbox{\vrule\kern1pt
       \vbox{\kern1pt#1\kern1pt}\kern1pt\vrule}\hrule}}
\def\cqfd{\hfill\boxit{\phantom{\i}}}
\def\cqfddd{\hfill\hbox{\cqfd\thinspace\cqfd\thinspace\cqfd}}

\newcount\numero \numero=1
\def\numeroter{{({\oldstyle\number\numero})}\ \advance\numero by1}

\newcount\refno 
\long\def\ref#1:#2<#3>{                                        
\global\advance\refno by1\par\noindent                              
\llap{[{\bf\number\refno}]\ }{#1} \pointir{\it #2} #3\goodbreak }

\newcount\refno 
\long\def\ref#1:#2<#3>{                                        
\global\advance\refno by1\par\noindent                              
\llap{[{\bf\number\refno}]\ }{#1} \pointir{\it #2} #3\goodbreak }

\def\citer#1(#2){[{\bf\number#1}\if#2\empty\relax\else,\ {#2}\fi]}

\def\boxit#1{\vbox{\hrule\hbox{\vrule\kern1pt
       \vbox{\kern1pt#1\kern1pt}\kern1pt\vrule}\hrule}}
\def\cqfd{\hfill\boxit{\phantom{\i}}}

\newbox\bibbox
\setbox\bibbox\vbox{\bigbreak
\centerline{{\pc BIBLIOGRAPHY}}

\ref{\pc BARTEL} (A) \& {\pc DOKCHITSER} (T):
Brauer relations in finite groups,
<J.\ Eur.\ Math.\ Soc. {\bf 17} (2015) 10, 2473--2512. >
\newcount\bartdok \global\bartdok=\refno


\ref{\pc DALAWAT} (C):
Serre's ``\thinspace formule de masse\thinspace'' in prime degree,
<Monats\-hefte Math.\ {\bf 166} (2012) 1, 73--92.
Cf.~arXiv\string:1004.2016v6.>
\newcount\monatshefte \global\monatshefte=\refno

\ref{\pc DALAWAT} (C):
Solvable primitive extensions,
<arXiv\string:1608.04673.>
\newcount\solprimp \global\solprimp=\refno

\ref{\pc DALAWAT} (C):
$\F_p$-representations over $p$-fields,
<arXiv\string:1608.04181.>
\newcount\irrepp \global\irrepp=\refno

\ref{\pc DALAWAT} (C):
Little galoisian modules,
<arXiv\string:1608.04182.>
\newcount\littlegalmods \global\littlegalmods=\refno

\ref{\pc DALAWAT} (C) \& {\pc LEE} (JJ):
Tame ramification and group cohomology,
<arXiv\string:1305.2580.>
\newcount\dalawatlee \global\dalawatlee=\refno

 \ref{\pc DEL \pc CORSO} (I) \& {\pc DVORNICICH} (R):
 The compositum of wild extensions of local fields of prime degree,
 <Monatsh.\ Math.\ {\bf 150} (2007) 4, 271--288.>
 \newcount\delcorso \global\delcorso=\refno

\ref{\pc DEL \pc CORSO} (I), {\pc DVORNICICH} (R) \& {\pc MONGE} (M):
On wild extensions of a $p$-adic field,
<J.\ Number Theory {\bf 174} (2017), 322--342.  Cf.~aXiv\string:1601.05939.> 
\newcount\deldvomonge \global\deldvomonge=\refno

\ref{\pc DIARRA} (B): 
Construction des extensions primitives d'un corps $p$-adique.  <Groupe
de travail d'analyse ultram\'etrique {\bf 9} (1981--82) 2, Expos\'e
24, 19 p.>
\newcount\diarra \global\diarra=\refno

\ref{\pc FONTAINE} (J-M):
Extensions finies galoisiennes des corps valu\'es complets \`a
valuation discr\`ete, <S\'eminaire Delange-Pisot-Poitou, Th{\'e}orie
des nombres {\bf 9} (1967--68) 1, Expos\'e 6, 21 p.>
\newcount\fontaine \global\fontaine=\refno

\ref{\pc HENNIART} (G):
Representations du groupe de Weil d'un corps local,
<sites.mathdoc.fr/PMO/feuilleter.php?id=PMO\_1979>
\newcount\henniart \global\henniart=\refno

\ref{\pc IWASAWA} (K):
On Galois groups of local fields.
<Trans.\ Amer.\ Math.\ Soc.\ {\bf 80} (1955), 448--469.>
\newcount\iwasawa \global\iwasawa=\refno


\ref{\pc KRASNER} (M):
Sur la primitivit\'e des corps\/ $\goth P$-adiques, <Mathematica,
Cluj {\bf 13} (1937) 4, 72--191. Cf.\
alpha.math.uga.edu/\string~lorenz/articles.html>
\newcount\krasner \global\krasner=\refno


\ref{\pc KRASNER} (M):
Le nombre de sur-corps primitifs d'un degr\'e donn\'e et le nombre de
sur-corps m{\'e}tagaloisiens d'un degr\'e donn\'e d'un corps de
nombres $\goth g$-adiques,
<C.\ R.\ Acad.\ Sc.\ {\bf 206} (1938) A, 876--877.>
\newcount\krasnercras \global\krasnercras=\refno


\ref{\pc KRASNER}  (M):
Nombre des extensions d'un degr\'e donn\'e d'un corps\/ $\goth P$-adique,
<Les tendances g\'eom.\ en algèbre et th\'eorie des nombres, CNRS,
Paris, 1966,  p.~143--169.> 
\newcount\krasnertendances \global\krasnertendances=\refno


\ref{\pc PATI} (M):
Extensions of degree $p^l$ of a $p$-adic field,
<Annali di Matematica Pura ed Applicata 2 June 2016, 1--21.
Cf.~arXiv\string:1511.02040.> 
\newcount\pati \global\pati=\refno

\ref{\pc SERRE} (J-P):
Corps locaux,
<Publications de l'Universit{\'e} de Nancago {\sevenrm VIII}, Hermann,
Paris, 1968, 245 p.>
\newcount\corpslocaux \global\corpslocaux=\refno


\ref{\pc WEIL} (A):
Exercices dyadiques, 
<Invent.\ Math.\ {\bf 27} (1974), 1--22.>
\newcount\weil \global\weil=\refno

} 

\centerline{\bf Wildly primitive extensions} 
\bigskip\bigskip 
\centerline{Chandan Singh Dalawat} 
\centerline{Harish-Chandra Research Institute, HBNI}
\centerline{Chhatnag Road, Jhunsi, Allahabad 211019, India} 
\centerline{\tt dalawat@gmail.com}

\bigskip\bigskip

{{\bf Abstract}.  A finite separable extension of a field is called
primitive if there are no intermediate extensions.  The most
interesting primitive extensions of a local field with finite residue
field are the wildly ramified ones, and our aim here is to parametrise
them in a canonical manner.

\footnote{}{{\it MSC2010~:} Primary 11S15, 11S37}
\footnote{}{{\it Keywords~:} Local fields, primitive extensions, wild ramification
}}

\bigskip\bigskip
\rightline{\it Meinem Doktorvater gewidmet}

\bigbreak
{\bf 1.  Introduction}
\bigskip

\numeroter  Let $p$ be a prime number and let $K$ be a $p$-field,
namely a local field with finite residue field of characteristic~$p$.
All extensions of $K$ appearing below are assumed to be {\it
separable} over $K$.  A finite extension $E$ of $K$ is called {\it
primitive\/} if $[E:K]>1$ and if the only extensions of $K$ in $E$ are
$K$ and $E$. It is easy to see that a tamely ramified extension of $K$
is primitive if and only if it is either unramified of prime degree or
totally ramified of prime degree $l\neq p$~; the latter are
parametrised by sections of the projection $K^\times\!/K^{\times
l}\to\Z/l\Z$ coming from the valuation on $K$.  Thus tamely ramified
primitive extensions are easy to classify.

We are interested in primitive $p$-extensions of the $p$-field $K$ ---
those whose degree is a power of~$p$.  If a primitive $p$-extension
$E$ is ramified, then clearly it is wildly and totally ramified over
$K$~; we say that $E$ is {\it wildly primitive\/} for short.  The only
primitive $p$-extension of $K$ which is not wildly primitive is the
unramified extension of degree~$p$ over $K$.

\numeroter   Extensions of degree~$p$ over $K$ are always
primitive (and wildly so if ramified).  They have been parametrised by
Del Corso and Dvornicich \citer\delcorso() if $K$ has
characteristic~$0$ and by the present author \citer\monatshefte() in
general.  More recently, Del Corso, Dvornicich and
Monge \citer\deldvomonge() have studied wildly primitive extensions of
$K$ of degree $p^n$ when $\char K=0$~; see also Pati \citer\pati()
when moreover $n$ is prime.  Our aim here is to parametrise all
primitive $p$-extensions of an arbitrary $p$-field $K$ by generalising
one of the main results of \citer\monatshefte() (the case $n=1$) and
its proof, and to compute their discriminants.  For some historical
remarks, see \S{\bf 9}.

\bigbreak
{\bf 2.  Notations and results}
\medskip

\numeroter Let $k$ be the residue field of $K$ and $q=\Card k$.  For
every $n>0$, put $e_n=p^n-1$, $K_n=K(\!\root e_n\of1)$ (where $\root
e_n\of1$ stands for an element of order $e_n$ in the multiplicative
group), and $L_n=K_n(\!\root e_n\of{K_n^\times})$, so that $K_n$ is
the unramified extension of $K$ of degree equal to the order $s_n$ of
$\bar q\in(\Z/e_n\Z)^\times$, and $L_n$ is the maximal abelian
extension of $K_n$ of exponent dividing $e_n$.  The ramification index
(resp.~the residual degree) of $L_n$ over $K$ is $e_n$
(resp.~$s_ne_n$).

Note that $L_n$ is a tamely ramified galoisian extension of $K$ of
group $G_n=\Gal(L_n|K)$~; it is split over $K$ in the sense that the
short exact sequence $\{1\}\to T\to G_n\to G_n/T\to\{1\}$ has a section,
where $T$ is the inertia subgroup of $G_n$.  If $K$ has
characteristic~$0$, then the $p$-torsion subgroup ${}_pL_n^\times$ of
$L_n^\times$ has order~$p$ (because $L_1$ contains $\!\root
p-1\of{-p}\;$).

View $\overline{L_n^\times}=L_n^\times\!/L_n^{\times p}$
(resp.~$\overline{L_n^+}=L_n^+\!/\wp(L_n^+)$, where $\wp(x)=x^p-x$) as
an $\F_p[G_n]$-module if $K$ has characteristic~$0$ (resp.~$p$).  Our
first main result is the following parametrisation of the set of
primitive $p$-extensions of $K$ of fixed degree~:

\numeroter  {\it The set of\/ $K$-isomorphism classes of primitive
extensions\/ $E$ of\/ $K$ of degree\/~$p^n$ is in canonical bijection
with the set of simple submodules\/ $D$ of the\/ $\F_p[G_n]$-module\/
$L_n^\times\!/L_n^{\times p}$, resp.~$L_n^+\!/\wp(L_n^+)$, of\/
degree\/~$n$, under the correspondence\/ $EL_n=L_n(\!\root
p\of{D})$, resp.~$EL_n=L_n(\wp^{-1}(D))$. }

In characteristic~$0$, this is a more precise version of the main
result of \citer\deldvomonge(3.2), as we specify the extension $L_n$
explicitly.  The proof (\S{\bf 6}) is a generalisation from the case
$n=1$ treated in \citer\monatshefte().  As there, we also determine
the structure of the filtered $\F_p[G_n]$-module
$L_n^\times\!/L_n^{\times p}$, resp.~$L_n^+\!/\wp(L_n^+)$, in \S{\bf
7}.

Later we shall define the {\it level\/} of a simple submodule $D$ of
$\overline{L_n^\times}$ or of $\overline{L_n^+}$ in terms of
the natural filtration on these $G_n$-modules.  Our second main result
relates the level of $D$ to the differental exponent of the
corresponding primitive $p$-extension $E$ of $K$ ($\oldstyle33$).
This allows us to compute --- in principle --- the number of primitive
extensions $E$ of $K$ of degree~$p^n$ with a given differental
exponent.

This parametrisation is illustrated in \S{\bf 8} in the simplest cases
of primitive quartic or octic extensions of dyadic fields ($p=2$,
$n=2$ or~$3$).

\numeroter  {\it Remarks}.  Let $\tilde K$ be a maximal galoisian
extension of $K$ containing the extensions $L_n$ $\oldstyle(3)$, and
let $M_n$ be the maximal abelian extension of $L_n$ in $\tilde K$ of
exponent $p$, so that the direct limit
$V=\lim\limits_{\longrightarrow}L_n$ is the maximal tamely ramified
extension of $K$ in $\tilde K$ and the direct limit
$B=\lim\limits_{\longrightarrow}M_n$ is the maximal abelian extension
of $V$ in $\tilde K$ of exponent~$p$, namely\/ $B=V(\!\root
p\of{V^\times})$, resp.\/ $B=V(\wp^{-1}(V))$.  Let $W$ be the
compositum of all wildly primitive extensions $E$ of $K$ in $\tilde
K$, or equivalently the compositum in $\tilde K$ of their galoisian
closures $\hat E$ over $K$.  It follows from ($\oldstyle4$) that if
$[E:K]=p^n$, then $\hat E\subset M_n$, and hence $W\subset B$.  It is
likely that $W=B$, just as the compositum of all degree-$p$ extensions
of $K$ in $\tilde K$ is $M_1$ \citer\monatshefte(Proposition~33).


\numeroter  Let $\Gamma=\Gal(V|K)$.  The structure of the
$\F_p[[\Gamma]]$-module $\Gal(B|V)$ has been determined
in \citer\littlegalmods($\oldstyle39$),
resp.~\citer\littlegalmods($\oldstyle42$), by showing that the dual
$\F_p[[\Gamma]]$-module $\overline{V^\times}=V^\times\!/V^{\times p}$,
resp.~$\overline{V^+}=V^+\!/\wp(V^+)$, is isomorphic to
$\F_p[[\Gamma]]^{[K:\Q_p]}\oplus\F_p$, resp.~$\F_p[[\Gamma]]^{(\N)}$.

\bigbreak
{\bf 3.  Solvable primitive $l$-extensions}
\medskip

\numeroter Let us recall from \citer\solprimp() a general algebraic
result which we need.  Fix a field $F$ and a maximal galoisian
extension $\tilde F$ of $F$.  All extensions of $F$ appearing below
are assumed to be contained in $\tilde F$.  A finite extension $E$ of
$F$ is called {\it solvable\/} if the group $G=\Gal(\hat E|K)$ is
solvable, where $\hat E$ is the galoisian closure of $E$ over $F$.
Galois proved that if $E$ is a solvable primitive ($\oldstyle1$)
extension of $F$, then $[E:K]=l^n$ for some prime $l$ and some $n>0$.

\numeroter Also fix the prime $l$ and the integer $n>0$.  Let $N$ be a
minimal normal subgroup of $G$.  As $[E:F]=l^n$, the $\F_l$-dimension
$N$ is~$n$.  The group $\Gal(\tilde F|F)$ acts on $N$ through its
quotient $G/N$~; the resulting $\F_l$-representation $\rho$ of
$\Gal(\tilde F|F)$ is irreducible and its image is solvable.  The
extension $E$ of $F$ is uniquely determined (up to $F$-isomorphism) by
the pair consisting of $\rho$ and the extension $\hat E$ of $L$, where
$F\subset L\subset \hat E$ is such that $N=\Gal(\hat E|L)$.  More
precisely,

\numeroter {\it Let\/ $F$ be a field, $l$ a prime number, and\/ $n>0$
an integer.  The map sending a primitive solvable extension\/ $E$ of\/
$F$ of degree\/ $l^n$ to its galoisian closure\/ $\hat E$ over\/ $F$
establishes a bijection between the set of\/ $F$-isomorphism classes
of such\/ $E$ with the set of pairs\/ $(\rho,M)$ consisting of an
irreducible degree-$n$ $\F_l$-representation\/ $\rho$ of\/
$\Gal(\tilde F|F)$ with solvable image and an abelian extension\/ $M$
of exponent\/~$l$ and degree\/~$l^n$ of the fixed field\/
$L_\rho=\tilde F^{\Ker(\rho)}$ of\/ $\rho$ such that\/ $M$ is galoisian
over\/ $F$ and the resulting conjugation action of\/ $\Gal(L_\rho|F)$
on\/ $\Gal(M|L_\rho)$ is given by $\rho$.} \cqfd

\numeroter In particular, given a primitive $l$-extension $E$ of
degree $l^n$ over $F$, there is a finite galoisian extension $L$ of
$F$, uniquely determined by $E$, such that $\hat E=EL$, $\Gal(\hat
E|L)$ is an $\F_l$-space of dimension~$n$, and the $\Gal(L|F)$-module
$\Gal(\hat E|L)$ is faithful and simple.

In what follows, we will take $F$ to be our $p$-field $K$ and $l$ to
be the prime~$p$.

\bigbreak
{\bf 4.  Sections and conjugates}
\medskip

\numeroter Another  purely algebraic ingredient we need is a lemma
used in \citer\deldvomonge(3.3), where its proof is attributed to
\citer\bartdok(6.1).  For the convenience of the reader, we briefly
reproduce it here.

\numeroter {\it Let\/ $l$ be a prime number.  Let\/ $G$ be a finite
group which has a normal subgroup\/ $A$ of order prime to\/ $l$ and
index a power of\/ $l$, and let\/ $C$ be a simple\/ $\F_l[G]$-module
of finite degree\/ $>1$.  Then\/ $H^1(G,C)=\{0\}$ and\/ $H^2(G,C)=\{0\}$. }

{\it Proof}.  The inflation-restriction sequence in this situation is
the exact sequence \citer\corpslocaux(Chapitre~VII, Proposition~4)
$$
\{0\}\to
H^1(G/A,C^A)\hfl{\hbox{Inf}}{}{10mm}
H^1(G,C)\hfl{\hbox{Res}}{}{10mm}
H^1(A,C).
$$
Since the orders of $A$ and $C$ are relatively prime, we have
$H^1(A,C)=\{0\}$.  The same reason, together with the hypotheses that
$G/A$ is an $l$-group and that $\dim_{\F_l}C>1$, implies that
$C^A=\{0\}$ and hence $H^1(G/A,C^A)=\{0\}$.  Therefore
$H^1(G,C)=\{0\}$.  This being so, the sequence
$$
\{0\}\to
H^2(G/A,C^A)\hfl{\hbox{Inf}}{}{10mm}
H^2(G,C)\hfl{\hbox{Res}}{}{10mm}
H^2(A,C).
$$
is exact \citer\corpslocaux(Chapitre~VII, Proposition~5), and a
similar argument leads to the conclusion $H^2(G,C)=\{0\}$.  \cqfd

\numeroter  As the group $H^2(G,C)$ classifies extensions of $G$ by the
$G$-{\it module\/} $C$, and as $H^1(G,C)$ classifies sections of the
neutral extension of $G$ by $C$ up to conjugation, we get the
following equivalent statement (in the multiplicative notation)~:
every extension\/ $\{1\}\to C\to\Gamma\to G\to\{1\}$ of\/ $G$ by\/ $C$
admits a section\/ $G\to\Gamma$, and any two sections are conjugate by
an element of\/ $\Gamma$.

\bigbreak
{\bf 5.  Irreducible $\F_p$-representations of $\Gal(\tilde K|K)$}
\medskip

\numeroter  Let's return to our local field $K$ and recall that our
aim is to parametrise the set of primitive $p$-extensions of $K$.  The
general algebraic result of \S{\bf 3} leads us to classify irreducible
degree-$n$ $\F_p$-representations $\rho$ of $\Gal(\tilde K|K)$, where
$\tilde K$ is a maximal galoisian extension of $K$.  Note that every
finite extension of $K$ is (separable by hypothesis ($\oldstyle1$)
and) solvable in the sense of ($\oldstyle7$), and the condition that
the image of $\rho$ be solvable is automatically satisfied.  In this
context, recall one of the main results from \citer\irrepp(),
employing the notation introduced in \S{\bf 2}.

\numeroter {\it Every irreducible\/ $\F_p$-representation of\/
  $\Gal(\tilde K|K)$ of degree\/~$n$ factors through the quotient\/
  $G_n=\Gal(L_n|K)$.}  \cqfd

\bigbreak
{\bf 6.  The proof of the parametrisation}
\medskip

\numeroter Before entering into the details, let us outline the
strategy of the proof of Theorem~($\oldstyle4$).  We will first show
that for every primitive extension $E$ of degree $p^n$ over $K$,

{\it a)\/ the extension\/ $EL_n$ of\/ $L_n$ is abelian of exponent\/
$p$ and degree\/ $p^n$,}

{\it b)\/ the extension\/ $EL_n$ of\/ $K$ is a galoisian, and}

{\it c)\/ the resulting\/ $\F_p[G_n]$-module\/ $C=\Gal(EL_n|L_n)$ is
simple.}

Recall ($\oldstyle3$) that the $p$-torsion subgroup ${}_pL_n^\times$
of $L_n^\times$ has order~$p$ if $K$ has characteristic~$0$.
Therefore, once we establish {\it a)}, there will be a unique
$n$-dimensional subspace $D$ of $L_n^\times\!/L_n^{\times p}$,
resp. of $L_n^+\!/\wp(L_n^+)$, such that $EL_n=L_n(\!\root p\of D)$ if
$\char K=0$ and $EL_n=L_n(\wp^{-1}(D))$ if $\char K=p$.  Once we
establish {\it b)}, we will know that $D$ is $G_n$-stable, and, once
we establish {\it c)}, we will know that the $\F_p[G_n]$-module $D$ is
simple, because there are canonical isomorphisms
$D=\Hom(C,{}_pL_n^\times)$, resp. $D=\Hom(C,\F_p)$.

We will then show that conversely,

{\it d) for every simple\/ $\F_p[G_n]$-submodule\/ $D$ of\/
$L_n^\times\!/L_n^{\times p}$ or of\/ $L_n^+\!/\wp(L_n^+)$, of
degree\/ $n$, there is a primitive extension\/ $E$ of\/ $K$
of degree\/ $p^n$, unique up to\/ $K$-isomorphism, such that\/
$EL_n=L_n(\!\root p\of D)$ or\/ $EL_n=L_n(\wp^{-1}(D))$ respectively.}

\numeroter Let's prove ($\oldstyle4$).  Let $E$ be a primitive
extension of $K$ of degree~$p^n$.  By the general algebraic theory of
\S{\bf 3} (applicable because $E$ is solvable in the sense of
($\oldstyle7$)), $E$ determines a finite galoisian extension $L$ of
$K$ such that the galoisian closure $\hat E$ of $E$ over $K$ is $\hat E=EL$,
the group $\Gal(\hat E|L)$ is an $\F_p$-space of dimension $n$, and it is
faithful and simple as a $\Gal(L|K)$-module.  By ($\oldstyle15$), we
have $L\subset L_n$~; in particular, $L$ is tamely ramified over $K$.
We claim that the extensions $\hat E$ and $L_n$ of $L$ are {\it linearly
disjoint\/}.

This is clear if $E$ is the unramified degree-$p$ extension of $K$ (in
which case $L=K$, $[\hat E:K]=p$, and $[L_1:K]=(p-1)^2$).  Otherwise, $\hat E$
is totally ramified of degree $p^n$ over $L$ whereas $L_n$ is tamely
ramified over $L$, and the claim follows.  Therefore $EL_n=\hat EL_n$ has
the properties {\it a), b)} and {\it c)} of ($\oldstyle16$).

\numeroter To establish the claim {\it d)} of $(\oldstyle16$), take
a simple submodule $D$ of $L_n^\times\!/L_n^{\times p}$ or of
$L_n^+\!/\wp(L_n^+)$ (according as $K$ has characteristic~$0$ or~$p$),
of dimension $n$ over $\F_p$.  Put $M=L_n(\!\root p\of D)$ or
$M=L_n(\wp^{-1}(D))$ respectively and put $C=\Gal(M|L_n)$.  Note that
$M$ is galoisian over $K$, and that the $\F_p[G_n]$-module $C$ is
simple because $C=\Hom(D,{}_pL_n^\times)$ or $C=\Hom(D,\F_p)$
respectively.

\numeroter  If $n=1$, so that the order of 
$C$ is $p$ and the order of $G_1$ is $(p-1)^2$, we have
$H^2(G_1,C)=\{0\}$ and $H^1(G_1,C)=\{0\}$.  Therefore the extension
$\Gal(M|K)$ of $G_1$ by $C$ splits, and any two sections are
conjugate.  In other words, there is a degree-$p$ extension $E$ of
$K$, unique up to $K$-isomorphism, such that $EL_1=M$, and we are
done.

\numeroter There is a similar argument when $n>1$.  Consider the maximal
 unramified extension $p$-extension $P$ of $K$ in $L_n$.  The subgroup
 $\Gal(L_n|P)$ of $G_n$ is obviously invariant under cojugation, of
 order prime to~$p$, and of index a power of $p$.  We can therefore
 apply ($\oldstyle13$) to our situation and conclude that the
 extension $\Gal(M|K)$ of $G_n$ by $C$ splits, and that any two
 sections are conjugate.  This means that there is a degree-$p^n$
 extension $E$ of $K$, unique up to $K$-isomorphism, such that
 $EL_n=M$.

\numeroter It remains to show that $E$ is primitive.  Indeed, suppose
there is an intermediate extension $K\subset F\subset E$.  Then $FL_n$
is galoisian over $K$ and $\Gal(M|FL_n)$ is a $G_n$-stable subspace of
$C$, therefore either $FL_n=L_n$ or $FL_n=M$, which implies that
either $F=K$ or $F=E$.  This completes the proof of Theorem~($\oldstyle4$)
following the strategy outlined in ($\oldstyle16$).  \cqfddd

\medbreak

\numeroter {\it Remarks}.  Another strategy for proving
($\oldstyle16$){\it d)} would be to consider $L=L_n^{\Ker(\rho)}$,
where $\rho$ is the action of $G_n$ on $C$, and to show directly that
there is an abelian extension $N$ of $L$ of exponent~$p$ and
degree~$p^n$, unique up to $L$-isomorphism, which is galoisian over
$K$ and such that $NL_n=M$.  If there is such an $N$, then
($\oldstyle9)$ implies the existence, uniqueness, and primitivity of
$E$.

\numeroter If $E$ is a wildly primitive extention of $K$ of
degree~$p^n$, $D$ its parameter, $M=EL_n$ (so that $M=L_n(\root p\of
D)$ or $M=L_n(\wp^{-1}(D))$ respectively), and $\rho$ the action of
$G_n=\Gal(L_n|K)$ on $D$, then the action of $G_n$ on the $\F_p$-space
$C=\Gal(M|L_n)$ is $\rho^*\otimes\omega$, where $\rho^*$ is the
contragradient (dual) of $\rho$ and $\omega:G_n\to\F_p^\times$ is the
character giving the action of $G_n$ on ${}_pL_n^\times$ (resp.~on
$\F_p$) if $K$ has characteristic~$0$ (resp.~$p$), just as in the case
$n=1$
\citer\monatshefte(Lemma~15). Note that the (wild)
ramification subgroup $\Gal(M|K)_1$ of $\Gal(M|K)$ is $\Gal(M|K)_1=C$.
It is also clear that the group $\Aut_K(E)$ is trivial unless $E$ is
cyclic over $K$ (of degree~$p$).

\numeroter We have parametrised primitive extensions $E$ of $K$ of
degree~$p^n$ by simple submodules $D$ of the $\F_p[G_n]$-module
$L_n^\times\!/L_n^{\times p}$ or $L_n^+\!/\wp(L_n^+)$ of degree~$n$.
We could equally well parametrise them by their galoisian closures
$\hat E$ as in ($\oldstyle9$), since the above proof characterises the
$\hat E$ which arise.  Indeed, a finite galoisian extension $F$ of $K$
is of the form $\hat E$ for some primitive extension $E$ of $K$ of
degree~$p^n$ if and only if, $F'$ being the maximal tamely ramified
extension of $K$ in $F$, two properties hold~:

{\it i)\/ the ramification subgroup\/ $H=\Gal(F|F')$ of\/ $G=\Gal(F|K)$ is
an\/ $\F_p$-space of dimension\/ $n$, and

\it ii)\/ the\/ $\F_p[G/H]$-module\/ $H$ is faithful and simple.}

Also, we need to look for $F'$ only among the subextensions of $L_n$.

\numeroter Clearly, a primitive {\it galoisian\/} extension
of any field is cyclic of prime degree.  We claim that if $E$ is a
primitive extension of the $p$-field $K$ whose {\it galoisian
closure\/} $\hat E$ is a {\it totally ramified $p$-extension\/} of
$K$, then $\hat E=E$ (and hence $E$ is ramified cyclic of degree~$p$
over $K$).  Indeed, the hypothesis on $\hat E$ implies that the
maximal tamely ramified extension of $K$ in $\hat E$ is $K$.  By the
preceding remark, $\Gal(\hat E|K)$ is a faithful simple
$\F_p[\Gal(K|K)]$-module.  It follows that $\hat E|K$ is cyclic of
degree~$p$, and hence $\hat E=E$.

\bigbreak
{\bf 7.  Little galoisian modules}
\medskip

\numeroter Our understanding of primitive extensions of $K$ cannot be
complete without working out the structure of the $\F_p[G_n]$-modules
$L_n^\times\!/L_n^{\times p}$ and $L_n^+\!/\wp(L_n^+)$, respectively
when $K$ has characteristic~$0$ and $p$.  This was determined by
Iwasawa \citer\iwasawa() (see also \citer\deldvomonge(4.4)) in
characteristic~$0$ and in \citer\littlegalmods() in general.  In
this \S, we recall these structure theorems and compute the
discriminant of a wildly primitive extension $E$ of degree $p^n$ over
$K$ in terms of its {\it parameter\/} $D$, the simple
$\F_p[G_n]$-submodule of $\overline{L_n^\times}$ or
$\overline{L_n^+}$, of degree~$n$, such that $EL_n=L_n(\!\root p\of
D)$ or $EL_n=L_n(\wp^{-1}(D))$, associated to $E$ by ($\oldstyle4$).

\smallbreak

\numeroter {\it Suppose that\/ $K$ is a finite extension of\/ $\Q_p$.
  The\/ $\F_p[G_n]$-module\/
  $\overline{L_n^\times}=L_n^\times\!/L_n^{\times p}$ is isomorphic
  to\/ ${}_pL^\times\oplus\F_p[G_n]^{[K:\Q_p]}\oplus\F_p$.}  \cqfd

\numeroter {\it Suppose that the $p$-field\/ $K$ has
  characteristic\/ $p$. The\/ $\F_p[G_n]$-module\/
  $\overline{L_n^+}=L^+\!/\wp(L^+)$ is isomorphic to\/
  $\F_p\oplus\F_p[G_n]^{(\N)}$.} \cqfd

\smallbreak

\numeroter What is important in both cases is the natural filtration
on the $\F_p[G_n]$-modules $\overline{L_n^\times}$ or
$\overline{L_n^+}$ which was studied in detail in
\citer\littlegalmods().  Consider the problem of computing the
discriminant of $E$ over $K$ in terms of $D$.  Put $M=EL_n$ and
$C=\Gal(M|L_n)$.  Note that the ramification filtration on $C$ has a
{\it unique ramification break\/} $\gamma(C)$ because the $G_n$-module
$C$ is simple and the ramification filtration is $G_n$-stable.

We denote by $\pgoth_{L_n}$ the unique maximal ideal of the ring of
integers of $L_n$, and we put $U_{L_n}^i=1+\pgoth_{L_n}^i$ for
every~$i>0$.


\numeroter For defining the {\it level\/} of the simple submodule $D$
of $\overline{L_n^\times}$ or $\overline{L_n^+}$ of degree~$n$,
suppose first that $K$ has characteristic~$0$ and let $e_{L_n}$ be the
ramification index of $L_n$ over $\Q_p$.

Notice that $e_{L_n}=c_n.(p-1)$ for some integer $c_n>0$.  The
filtration on $\overline{L_n^\times}$ is given by the images $\bar
U_{L_n}^i$ of $U_{L_n}^i$ for various~$i>0$.  Put $\bar
U_{L_n}^0=\overline{L_n^\times}$ by convention.  As $D$ is a {\it
simple\/} submodule of $\overline{L_n^\times}$, there is a {\it
unique\/} $i\in\N$ such that $D\subset\bar U_{L_n}^i$ but $D\cap\bar
U_{L_n}^{i+1}=\{1\}$, because the filtration is $G_n$-stable.  We
define the {\it level\/} of $D$ to be $\delta(D)=pc_n-i$.

We have $\delta(D)\in[0,pc_n]$, and if $\delta(D)\equiv0\pmod p$, then
$n=1$ and either $\delta(D)=0$, $D=\bar U_{L_1}^{pc_1}$ (which is
$G_1$-isomorphic to ${}_pL_1^\times$), and $E$ is  unramified of
degree~$p$ over $K$, or $\delta(D)=pc_1$~; the latter lines $D$ and
the corresponding extensions $M$ of $L_1$ and $E$ of $K$ are said to
be {\it tr\`es ramifi\'ees}.

\numeroter Now suppose that $K$ has characteristic $p$.  The
filtration on $\overline{L_n^+}$ is given by the images
$\overline{\pgoth_{L_n}^i}$ ($i\in\Z$).  As before, and for the same
reason, since $D$ is a {\it simple\/} submodule of $\overline{L_n^+}$,
there is a {\it unique\/} $i\in\Z$ such that
$D\subset\overline{\pgoth_{L_n}^i}$ but
$D\cap\overline{\pgoth_{L_n}^{i+1}}=\{0\}$.  We define the {\it
level\/} of $D$ to be $\delta(D)=-i$.

We have $\delta(D)\in\N$, and if $\delta(D)\equiv0\pmod p$, then
$n=1$, $\delta(D)=0$, $D=\overline{\pgoth_{L_1}^0}$ (which is
$G_1$-isomorphic to $\F_p$), and $E$ is unramified of degree~$p$ over
$K$.  There is no analogue of {\it tr\`es ramifi\'ees} lines or
extensions.

\numeroter For finite extensions $L$ of $K$ and $M$ of $L$, denote
the {\it differental exponent\/} (resp. ramification index) of $M|L$
by $d_{M|L}$ (resp.~$e_{M|L}$), and recall that
$d_{M|K}=d_{M|L}+d_{L|K}e_{M|L}$ \citer\corpslocaux(Chapitre~III,
Proposition~8).  If $M|L$ is galoisian of group $G$, then
$d_{M|L}=\sum_{i\in\N}(g_i-1)$, where $g_i$ is the order of the higher
ramification subgroup $G_i\subset G$ (in the lower
numbering) \citer\corpslocaux(Chapitre~IV, Proposition~4).  If $L$ is
tame over $K$, then
$d_{L|K}=e_{L|K}-1$ \citer\corpslocaux(Chapitre~III, Proposition~13).
So a good measure of the {\it wildness\/} of $L$ over $K$ in general
is $\varepsilon(L)=d_{L|K}-(e_{L|K}-1)$, and a good name for the
invariant $\varepsilon(L)$ would be the {\it differental excess\/} of
$L|K$~; it was used by Serre in his mass formula for totally ramified
extensions of $K$ of fixed degree.

\numeroter {\it Let\/ $E$ be a wildly primitive extension of\/ $K$,
$p^n=[E:K]$ its degree, $\varepsilon(E)$ its differental excess\/
  $\oldstyle(32)$, and\/ $\gamma(C)$ the unique ramification break\/
  $\oldstyle(29)$ of\/ $C=\Gal(M|L_n)$, where\/ $M=EL_n$.  Let\/ $D$
  be the simple\/ $\F_p[G_n]$-submodule of\/ $\overline{L_n^\times}$
  or\/ $\overline{L_n^+}$, of degree\/~$n$, such that\/ $M=L_n(\!\root
  p\of D)$ or\/ $M=L_n(\wp^{-1}(D))$, and\/ $\delta(D)$ its level as
  in\/ $\oldstyle(30)$ or\/ $\oldstyle(31)$.  We have\/
  $\gamma(C)=\delta(D)=\varepsilon(E)$.}

{\it Proof}.  See \citer\monatshefte(34) for the case $n=1$~; the same
proof works for $n>0$.  The equality $\gamma(C)=\delta(D)$ follows
from a certain orthogonality relation recalled there (where the level
of $D$ was defined to be $-\delta(D)$).  Apply ($\oldstyle32$) to
get\par
\hbox{\kern-1.5cm
\vbox{
$$
\diagram{
E&\hfl{e_n}{}{10mm}&M&\quad\qquad&E&\hfl{e_n-1}{}{10mm}&M\cr
\ufl{p^n}{}{5mm}&&\ufl{}{p^n}{5mm}&\quad\qquad&
 \ufl{d_{E|K}}{}{5mm}&&\ufl{}{(1+\gamma(C))(p^n-1)}{5mm}\cr
K&\hfl{}{e_n}{10mm}&L_n&\quad\qquad&K&\hfl{}{e_n-1}{10mm}&L_n,\cr}
$$}}
\noindent where the numbers along the arrows in the first
(resp.~second) square indicate ramification indices (resp.~differental
exponents) of the corresponding extension.  Compute $d_{M|K}$ along
the two paths from $K$ to $M$ and compare to get
$$
(1+\gamma(C))(p^n-1)+(e_n-1)p^n=(e_n-1)+d_{E|K}e_n,
$$
and  recall that $e_n=p^n-1$, to conclude that
$\gamma(C)=\varepsilon(E)$.  \cqfd

\bigbreak
{\bf 8.  Some quartic and octic examples}
\medskip

\numeroter  Taking $p=2$ and $n=2$, we will briefly indicate how  to recover
primitive quartic extensions $E$ of dyadic fields $K$ which were
studied by Weil \citer\weil().  We will say that $E$ is an ${\goth
A}_4$-quartic (resp.~${\goth S}_4$-quartic) if the group $\Gal(\hat
E|K)$ of its galoisian closure $\hat E$ is isomorphic to ${\goth A}_4$
(resp.~${\goth S}_4$).

\numeroter First let $K=\Q_2$.  In view of
($\oldstyle4$), we should look for $G_2$-stable irreducible
$\F_2$-planes in $\overline{L_2^\times}=L_2^\times\!/L_2^{\times 2}$.
(Recall that $G_2=\Gal(L_2|K)$ and $L_2$ is the maximal abelian
extension of exponent~$e_2=3$ of the unramified quadratic extension
$K_2=K(\root3\of1)$ of $K$).

But first let us classify degree-$2$ $\F_2$-representations of $G_2$.
Let $l_2$ be the residue field of $L_2$, let $T\subset G_2$ be the
inertia subgroup, and let $\theta:T\to l_2^\times$ be the canonical
character.  Then, in the notation of \citer\irrepp($\oldstyle19$), we
have

\numeroter {\it The only irreducible degree\/-$2$
$\F_2$-representations of\/ $G_2$ are\/
$\pi_{\overline{\bar1,\root3\of1}}$ and\/
$\pi_{\overline{\bar\theta,1}}$.} \cqfd

\numeroter  Concretely, $L_2$ contains the unramified cubic extension
$L=K(\root7\of1)$ of $K$, and $\pi_{\overline{\bar1,\root3\of1}}$ is
the unique irreducible degree\/-$2$ $\F_2$-representations of\/
$\Gal(L|K)=\Z/3\Z$.  Similarly, $L_2$ contains the
unique \citer\dalawatlee(8.1) ${\goth S}_3$-extension
$L'=K(\!\root3\of1,\root3\of2)$ of $K$, and
$\pi_{\overline{\bar\theta,1}}$ is the unique irreducible degree\/-$2$
$\F_2$-representations of\/ $\Gal(L'|K)={\goth S}_3$.  So it suffices
to work separately over  $L$ and $L'$  instead of
$L_2=LL'$.

\numeroter  First take $L=K(\!\root7\of1)$ and let $G=\Gal(L|K)$.  It
follows from \citer\littlegalmods($\oldstyle21$) that the
$\F_2[G]$-module $\bar U_L^1$ is isomorphic to
${}_2L^\times\oplus\F_2[G]$, which contains a unique $G$-stable
$\F_2$-plane $D$. By $\oldstyle(9)$, there is a unique primitive
quartic extension $E$ of $K$ whose galoisian closure is $\hat
E=L(\!\root2\of D)$~; the group $\Gal(\hat E|K)$ is isomorphic to
${\goth A}_4$, and $\hat E$ is the unique ${\goth A}_4$-extension of
$K$.

\numeroter  Now let  ${L'}=K(\!\root3\of1,\root3\of2)$ and $G=\Gal({L'}|K)$.
It follows as before that the $\F_2[G]$-module $\bar U_{L'}^1$ is
isomorphic to ${}_2{L'}^\times\oplus\F_2[G]$, and a finer analysis as
in \citer\littlegalmods($\oldstyle11$) shows that there is a unique
$G$-stable irreducible $\F_2$-plane $D\subset\bar U_{L'}^5$ (such that
$D\cap\bar U_{L'}^6=\{1\}$), and two $G$-stable irreducible $\F_2$-planes
$D\subset\bar U_{L'}^1$ such that $D\cap\bar U_{L'}^2=\{1\}$.
Corresponding to each $D$, we get a primitive quartic extensions $E$
of $K$ whose galoisian closure is $\hat E={L'}(\!\root2\of D)$~; the
group $\Gal(\hat E|K)$ is isomorphic to ${\goth S}_4$ in each case,
and these three are the only ${\goth S}_4$-extensions of $K$.

\numeroter  The differental exponents  of these
extensions can be computed as in ($\oldstyle33$).  We leave for the
readers (or their indefatigable computers) the pleasure of doing so
and of finding equations defining them.  See \citer\henniart(p.~111)
for the details.  

\numeroter  Similar computations can be made for
$K=\F_2\series{\varpi}$, by working over the unramified cubic
extension $K(\!\root7\of1)$ (resp.~the unique ${\goth S}_3$-extension
$K(\!\root3\of1,\root3\of\varpi)$).  Now there are infinitely many
${\goth A}_4$-quartic (resp.~${\goth S}_4$-quartic) extensions, but
only finitely many with bounded differental exponent.  The ${\goth
A}_4$-quartics have to be counted carefully (as
in \citer\irrepp($\oldstyle25$)), because the corresponding
$\F_2$-representation is not absolutely irreducible~: over the
quadratic extension $\F_2(\!\root3\of1)$ of $\F_2$, it splits into the
direct sum of the two cubic characters
$\Gal(K(\!\root7\of1)|K)\to\F_2(\!\root3\of1)^\times$.  These two
characters are interchanged by the generator of the group
$\Gal(\F_2(\!\root3\of1)|\F_2)$ of order~$2$.

\numeroter Finally allow $K$ to be any $2$-field and let $q$ be the
cardinal of its residue field.  If $q\equiv-1\pmod3$, the theory is
completely similar to the cases $q=2$ discussed above.  If
$q\equiv1\pmod3$ (so that $K$ contains $\root3\of1$), then the group
$G_2$ is commutative of exponent~$3$ (and order $3^2$).  So apart from
the unramified cubic extension, we have to deal with the three
ramified cubic extensions (all three cyclic) of $K$.  Each gives a
certain (finite) number of ${\goth A}_4$-quartic extensions of $K$ of
bounded differental exponent.  There are no ${\goth S}_4$-quartic
extensions because $K$ has no ${\goth S}_3$-extensions.

\numeroter Taking $K=\Q_2$ and $n=3$, one can recover
the list of primitive octic extensions of $\Q_2$ and their differental
exponents to be found in \citer\krasnercras().  We have $e_3=2^3-1=7$,
$K_3=K(\!\root7\of1)$, $L_3=K_2(\root7\of\xi,\root7\of2)$ (where $\xi$
is a generator of the multiplicative group $k_3^\times$ of the residue
field $k_3$ of $K_3$), and $G_3=\Gal(L_3|K)$.  One has to determine
the irreducible degree-$3$ $\F_2$-representations $\pi$ of $G_3$, the
copies of each $\pi$ in $L_3^\times\!/L_3^{\times 2}$, and the level
($\oldstyle30$) of each copy.  The copies correspond to primitive
octic extensions of $K$ by ($\oldstyle4$) and the levels are related
to their differental exponents by ($\oldstyle33$).  The same
computation works for $\F_2\series{\varpi}$ upon replacing
$L_3^\times\!/L_3^{\times 2}$ with $L_3^+/\wp(L_3^+)$.  We omit the
details.

%
%
%
%

\bigbreak
{\bf 9.  Historical note}
\medskip

\numeroter Let us finish by listing a few papers related to the theme of
primitivity not already mentioned in the Introduction.  This account is
far from being a history of the subject.  The concept of primitivity
for subgroups of the symmetric group goes back to the {\it Second
m{\'e}moire\/} (1830) of Galois~; it was clarified by Jordan in his
thesis (1860) and in his {\it Trait\'e} (1870).

\numeroter In a long series of papers beginning with \citer\krasner()
and culminating in \citer\krasnertendances(), and in several notes in
the {\it Compte rendus}, Krasner appears to have been the first to
study wildly primitive extensions, initially over a finite extension
$K$ of $\Q_p$ and later also over $p$-fields of characteristic $p$.
He introduces the notion of hypergroups and extends ramification
theory to finite extensions $F$ of $K$ which are not assumed to be
galoisian over $K$, and, according to Arf's review of \citer\krasner()
in the {\it Zentralblatt\/} {\bf 18} (p.~202), gives a necessary and
sufficient condition for $F$ to be primitive.  To illustrate his
theory, he computes Eisenstein polynomials defining the sixteen
primitive octic extensions $E$ of $\Q_2$ \citer\krasnercras(), along
with information from which the differental exponent of $E$ can be
deduced. His method was taken up and generalised by Diarra
\citer\diarra() to $p$-fields of characteristic~$p$.  I haven't
succeeded in penetrating their work.

\numeroter  We have seen ($\oldstyle28$) that for a wildly
primitive extension $E$ of a $p$-field $K$ with galoisian closure
$\hat E$, the wild ramification subgroup of $\Gal(\hat E|K)$ has a
unique ramification break~; in particular, it is commutative of
exponent~$p$.  In \citer\fontaine(), Fontaine studies Eisenstein
polynomials over a $p$-field which define an abelian extension of
exponent~$p$ with a unique ramification break.

\numeroter In \citer\weil(), Weil  studies ${\goth A}_4$- 
and ${\goth S}_4$-extensions of dyadic fields $K$.  As we have seen,
they are the same as galoisian closures of primitive quartic
extensions, so in principle their enumeration follows from the work of
Krasner and Diarra.  Weil's results can be viewed as the case $p=2$,
$n=2$ of the foregoing.

\bigbreak
{\bf 10.  Acknowledgements}. Work on this project of parametrising
wildly primitive extensions of $p$-fields, generalising from the case
of degree-$p$ extensions treated earlier \citer\monatshefte(), was
begun when the author was enjoying the hospitality of the Research
Institute for Mathematical Sciences, Kyoto, and he would like to thank
Akio Tamagawa and Kyoko Price for making the stay so fruitful. This
Note completes the
sequence \citer\solprimp()--\citer\littlegalmods()~; all four papers
have been influenced by \citer\deldvomonge() at various places.  I am
extremely grateful to Dino Lorenzini for making \citer\krasner()
available on his website in 2011.

\bigbreak
\unvbox\bibbox 

\vfill\eject
\bye